\documentclass[12pt]{article}
\usepackage{graphicx}
\usepackage{amsmath}\usepackage{mathrsfs,amsfonts}
\usepackage{amsfonts}
\usepackage{amsmath,amssymb,amsfonts}
\usepackage{amssymb}
\usepackage{color}

\usepackage{dsfont}
\usepackage{paralist}

\usepackage{bm}
\usepackage{epsfig}

\definecolor{c}{rgb}{0.9,0.3,0.1}
\definecolor{b}{rgb}{0.1,0.3,0.9}

\newtheorem{remark}{Remark}[section]
\newtheorem{lemma}{Lemma}[section]
\newtheorem{theorem}{Theorem}[section]

\newtheorem{hypothesis}{Hypothesis}[section]
\renewcommand{\theequation}{\arabic{section}.\arabic{equation}}

\def\de{{\delta}}
\def\ep{{\epsilon}}
\def\ga{{\gamma}}
\def\si{{\sigma}}\def\th{{\theta}}\def\ze{{\zeta}}\def\va{{\varepsilon}}

\def\<{\left<}\def\>{\right>}\def\({\left(}\def\){\right)}

\newfam\msbmfam\font\tenmsbm=msbm10\textfont
\msbmfam=\tenmsbm\font\sevenmsbm=msbm7
\scriptfont\msbmfam=\sevenmsbm\def\bb#1{{\fam\msbmfam #1}}

\def\EE{\bb E}

\newcommand{\cE}{{\mathcal{E}}}
\newcommand{\cG}{{\mathcal{G}}}
\newcommand{\cH}{{\mathcal{H}}}
\newcommand{\cU}{{\mathcal{U}}}
\newcommand{\cX}{{\mathcal{X}}}

\numberwithin{equation}{section}

\begin{document}
\title{
A stochastic maximum principle for partially observed stochastic control systems with delay \footnote{SZ's research is supported in   by 
the Fundamental Research Funds for the Central Universities (Grant No. 2020ZDPYMS01)
and Nature Science Foundation of China (Grant No. 11501129). 
 LX's research is supported partially by RGC (Grant No. 15213218 and 15215319).
 JX's research is supported by Southern University of Science and Technology Start up fund Y01286120 and National Natural Science Foundation of China (grant nos. 61873325, 11831010). Corresponding author: Shuaiqi Zhang}}
\author{  {Shuaiqi Zhang}\\ \small
 School
of Mathematics, China  University of Mining and Technology, 
\\ \small
Xuzhou, Jiangsu, 221116, China\\
\\
 {Xun Li}\\ \small
Department  of Applied Mathematics, The Hong Kong Polytechnic University,
\\ \small
Hong Kong, China,\\
\\
{Jie  Xiong
}\\  \small Department of Mathematics, Southern University of Science and Technology,
\\
 \small Shenzhen, China\\
    }
\maketitle
\date

\begin{abstract}
This paper deals with  partially-observed optimal control problems for the state governed
by stochastic differential equation with delay. 
We develop a stochastic maximum principle for this kind of optimal control problems using a variational method and a filtering technique. 
Also, we establish a sufficient condition without assumption of the concavity.
Two examples shed light on the theoretical results are established in the paper. In particular, in the example of an optimal investment problem with delay, 
its numerical simulation shows the effect of delay via a discretization technique for forward-backward stochastic differential equations (FBSDEs) with delay and anticipate terms.  


{\bf AMS 2000 Subject Classification:}
\end{abstract}

\section{Introduction}\label{intro}
\setcounter{equation}{0}
\renewcommand{\theequation}{1.\arabic{equation}}

This paper is concerned with partially-observed optimal control problems for stochastic differential equations  with delay.
This kind of optimal control problems has
a variety of important applications in many fields such as science, engineering, economics and finance.
The past decades witnessed tremendous interests and research efforts in dealing with partially observable problems using variational methods as well as their applications;
see e.g., Bensoussan \cite{ba2} studied stochastic control of partially observable systems,
Zhang \cite{zq} formulated a class of controlled problem on partially observed diffusions with correlated noise, 
Zhou \cite{zx} extablished the necessary conditions of optimal controls for stochastic partial differential equations, 
Tang \cite{ts} analyzed the maximum principle for partially observed optimal control of stochastic differential equations, Wang and Wu \cite{WW} derived Kalman-Bucy filtering equations of forward and backward stochastic equations (FBSDEs) and considered applications to recursive optimal control problems,
Wang $et\; al$ \cite{wzz} further developed stochastic maximum principle for mean-field type optimal control with partial information.
In these papers, however, the information filtration does not depend on the control itself.
Other works with coupled information filtration and control include Huang $et\; al$ \cite{hwx} considering partial information control problems for backward stochastic differential equations (BSDEs), Wang $et\; al$ \cite{wwx} and \cite{wwx2} studying linear-quadratic  control problems based on FBSDEs. In addition, Zhang $et\; al$ \cite{zhang2017} investigated FBSDEs with jumps and regime switching. 
However, all these research works do not involve the delay term.
In summary, this is still a fascinating yet challenging research area when we consider the stochastic system with delay.

On the other hand, taking into account of the fact that it may take some time before new information affects the state, and the control may based on both present and past values of the state, a delay term should be considered in the controlled systems. 
Such problems arise in many applications in the fields of biology, mathematical finance,  physics, engineering and so on. For instance,
there is a time lag between the wealth process and regulation or the after effect of investment.
 As for the optimal control problem with delay,  {\O}ksendal $et\; al$ \cite{OksendalSulemZhang2011} investigated
the  controlled system with delay and jumps in the state equation. They established sufficient and necessary conditions for stochastic maximum principles. Chen and Wu \cite{chenwu2010}  studied the stochastic maximum principle for a system involving both delays in the state variable and the control variable.
Huang $et\; al$ \cite{HLS2012}
 discussed the optimal control for
FBSDEs with delay in forward equation and anticipated term in backward equation.
However, these papers deal with the problem within the framework of full information.

Recently, there has been increasing interest in studying partially-observed control problems with delay as well as addressing their applications. 
However, the topic of delay-based optimal partially-observed control problems is still a relatively under-explored research field, and therefore many fundamental questions remain open and new methodologies need to be developed. A recent paper we should mention is Li $et\; al$ \cite{LWW} in which the authors studied the linear-quadratic optimal control for time-delay stochastic system with recursive utility under partial information. 
In this paper, we focus on studying partially-observed control problems for nonlinear state equation with delay, where information filtration, generated by the observation process, depends on the control in a coupled manner.
We formulate and tackle the problem within the framework of filtering theory corresponding to the specific model of partially-observed differential equation with delay. Compared with the aforementioned work in the setup of full information, the control is adapted to the filtration generated by the observation process, rather than the natural filtration generated by  Brownian motion and poisson process, resulting in substantial difference in the analysis.
Last but not least, the forward stochastic differential equation with delay describes the state process while the adjoint equation is a BSDE with anticipated term. Since the state equation involves delay term (which implies path-dependency), the system is not Markovian. This means that the classical four step scheme no longer works and its solution cannot be derived explicitly and directly.  Therefore, we adopt Euler scheme to study its numerical solution.

The rest of this paper is organized as follows. We formulate the partial information stochastic control problem in Section 2. In Section 3, a necessary  condition  for optimal control of   stochastic systems with delay under partial information is established. Section 4 is concerned with the sufficient condition.
Two examples are provided for illustrative purpose in Section 5. 
Also, simulation for our model using the numerical scheme  developed in \cite{zhang2019} shows that the delay term causes a decrease in  the value function.
Throughout this paper, we will use $K$ to denote a constant whose value can be changed from place to place.

\section{The model}
\setcounter{equation}{0}
\renewcommand{\theequation}{2.\arabic{equation}}

Let $\left(\Omega,   \mathcal{F} ,  \mathcal{F}_t,  \mathrm P  \right) $
 be a filtered  probability space on which   a real-valued standard
Brownian motion   $B(\cdot)$ is defined.
Suppose that $N(d\zeta, dt)$ is a Poisson random measure (independent of $B$) defined on $\mathrm R_+\times\cE$ where $(E,\cE,\nu)$ is a measure space. Now we define the compensated Poisson random measure as
 \[\widetilde{N} (d\zeta, dt)=N(d\zeta,dt)-\nu (d\zeta) dt, \]
where the intensity measure $\nu (d \ze )$  serves the purpose of
compensating.

Consider the stochastic differential equation with delay
\begin{eqnarray}\label{state}
\left\{\begin{array}{rcl}
dx (t) &=& b(t, x (t), x(t-\delta), v(t),v(t-\delta))dt \\
 &&+ \sigma(t, x (t), x(t-\delta), v(t),v(t-\delta) )dB(t)\\
&& + \theta(t, x (t), x(t-\delta), v(t),v(t-\delta) )dW^v(t)\\
&& + \!\int_{\mathcal{E}}\!\gamma (t, x (t-), x(t\!-\!\delta), v(t-),v(t\!-\!\delta), \zeta )\widetilde{N}(d\zeta, dt), \\
&& 
\hfill t\in[0,T],     \\
x (t) &=& \xi(t), \quad v(t)=\eta(t),\quad\quad  t\in[-\delta,0],
\end{array}\right. \nonumber \\
\end{eqnarray}
where $v(t)$ is a control process taking values in a convex set $U \subseteq \mathrm R$,
and $W^v$, taking values in $\mathrm R$, is a stochastic process depending on the control process $v(t)$, and $\delta$ is the time delay.
In the above equation,
 $b,\ \sigma,\ \th:[0,T]\times \mathrm R \times \mathrm R \times U\times U \rightarrow   \mathrm R  $,
 and
 $\gamma:[0,T]\times\mathrm R   \times \mathrm R  \times U\times U \times \cE \rightarrow \mathrm R     $.

Suppose that the observation process $Y(t)$ is a Brownian motion which satisfies the following equation:

\begin{equation}\label{observation}
\left\{\begin{array}{rcl}
dY(t) & \!\!=\!\! & h(t, x (t))dt+dW^v(t),\\
Y(0)  & \!\!=\!\! & 0,
\end{array}
\right.
\end{equation}
where $h: [0,T]\times \mathrm R     \rightarrow  \mathrm R       $ is a given continuous mapping.

Substituting  (\ref{observation}) into (\ref{state}) yields
\begin{eqnarray}\label{subs-state}
\left\{\begin{array}{rcl}
dx (t) &=&  \tilde{b}(t, x (t), x(t-\delta), v(t),v(t-\delta))dt \\
&&+\sigma(t, x (t), x(t-\delta), v(t),v(t-\delta))dB(t) \\
&&+ \theta(t, x (t), x(t-\delta), v(t),v(t-\delta) )d Y(t) \\
&& +\!\int_{\mathcal{E}}\!\!\gamma(t, x (t-\!), x(t\!-\!\delta), v(t-),v(t\!-\!\delta), \zeta )\widetilde{N}(d\zeta, \!dt), \\
x (t) &=& \xi(t), \quad v(t)=\eta(t), \quad  t\in[-\delta,0],
\end{array}\right. \nonumber \\
\end{eqnarray}
where $\tilde{b}=b-\theta h$. Let
\begin{eqnarray}\label{subfil}
\cG_t= \bm  \si\(Y(s):\;0\le s\le t\).
\end{eqnarray}
Denote  the set of admissible controls by $\mathcal{U}_{ad}$ which consists of all
 $\mathcal{G}_t$-adapted $U$-valued processes s.t. $\mathbb{E}\int^T_0v(t)^8dt<\infty$
 and $v(t)=\eta(t), \  t\in[-\delta,0]$.

Note that we restrict our study to the scalar case for simplicity of notation only. The results can be modified to suit the vector case.

\begin{hypothesis}\label{HP1}
Fix $t$, the functions $b$, $\sigma$, $\theta$,    and $h$  are of bounded continuous
derivatives up to second order with respect to other variables;
 $\gamma $, as an $L^2(\cE,\nu)$-valued function, is twice Fr\'echet differentiable with respect to other variables, and its
    Fr\'echet derivatives are  bounded.
Also, we assume that
\begin{itemize}
\item  there exists a constant $K$ such that for all
$t_i,  x_i,x'_i, v_i, v'_i$, $i=1,2$,
\begin{eqnarray*}
\begin{array}{rl}
& |\phi(t_1, x_1,x_1',v_1,v_1')- \phi(t_2, x_2,x_2', v_2,v_2')|_\cX \quad \quad
 \quad \nonumber\\
  \leq  & K (|t_1-t_2|^{1/2}+ |x_1-x_2 |  +|x_1'- x_2'  |+|v_1 - v_2   |+ |v_1'- v_2' | ) ,
\end{array}
\end{eqnarray*}
 \item   \[  |\phi(t,0,0,0,0)   |_\cX\le K, \;\forall t\in[0,T], \]
 \end{itemize}
where  $\cX=\mathrm R$ for $ \phi=b,\ \sigma,\ \theta$ and $\cX=L^2(\cE,\nu)$ for $\phi=\ga$.

Moreover, there is a constant $ C$ such that
$ | h(t,x) | \leq C  $ for any $(t,x)\in  [0,1]\times \mathrm R  $.
\end{hypothesis}

 Let $ \bar{\mathrm P}$ be a probability measure defined as  $d \bar{\mathrm P}=Z^{v}(t) d \mathrm P$, where $Z^{v}(t)$ is given by
\begin{equation}\label{Z}
\left\{\begin{array}{rcl}
 dZ^{v}(t)&=&Z^{v}(t)h (t, x (t))dY(t),\\
Z^{v}(0)&=& 1.
\end{array}
\right.
\end{equation}

From Girsanov's theorem and (\ref{observation}), $(B(t),W^v(t))$ is a 2-dimensional  Brownian motion defined on
 $\left(\Omega,   \mathcal{F} ,  \mathcal{F}_t, \bar {\mathrm P}  \right) $.

The performance functional is assumed to have the form
\begin{eqnarray}\label{JP}
\begin{array}{l}
J(v(\cdot)) = \bar {\mathbb{E}}\Big[\int_{0}^{T}  \ell\left(t,
               x (t), v(t), v(t-\delta) \right)dt+   \varphi(x(T)) \Big],\\
\end{array}  
\end{eqnarray}
where the notation $\bar{\EE}$ stands for the expectation with respect to the probability
measure $\bar {\mathbb  {P}}$.

\begin{hypothesis}\label{HP2}
The functions $\ell$ and $\varphi $  are of linear growth: there exists a constant such that $\forall (t,x,v,v')$,
\[
|\ell(t,x,v,v')|\!+\!|\phi(x)|\le K(1+|x|+|v|+|v'|).
\]

By the definition of the admissible control and Lemma \ref{lemma1} below, it is easy to verify that for any $v\in\cU_{ad}$,
\[
\bar {\mathbb E} \bigg[\int_0^T |\ell \(t, x(t),  v(t), v(t-\delta) \)| dt+| \varphi \(x(T)\) | \bigg]  <\infty.
\]
\end{hypothesis}
\hfill $\Box$

Note that, by Girsanov's theorem, the performance functional (\ref{JP}) can be rewritten as
\begin{eqnarray}\label{JZ}
J(v(\cdot))&=& \EE \Big[\int_{0}^{T}Z^{v}(t)\ \ell\left(t,
               x (t), v(t), v(t-\delta) \right)dt\nonumber\\&&+Z^{v}(T)\ \varphi(x(T)) \Big].
\end{eqnarray}

 The  optimal control problem is to find a control $ v  $ to maximize  the  functional
(\ref{JP}) over $v(\cdot)\in \cU_{ad}$  subject to
(\ref{subs-state}) and (\ref{Z}),  i.e.,
\begin{equation}
 J( u)=\max_{v(\cdot)\in  \cU_{ad} } J (v(\cdot)).
 \nonumber\\
\end{equation}

\begin{remark}\label{RJ}
Usually, the cost functional $J(\cdot)$ should depend on $v(\cdot) $. In some cases, for instance, an insurance company adopts a reinsurance strategy to reduce its exposure to risk,  there is a delay in the reaction of the strategy. Therefore, the cost functional is related to $v(t-\delta)$. This consideration leads to the form (\ref{JZ}) for the cost functional.
\end{remark}

 Finally, we proceed to presenting our main results. We first fix an admissible control $u$ and denote \[b(t)=b(t,x(t),x(t-\de),u(t),u(t-\de))\]
and
\[
b_x(t)=b_x(t,x(t),x(t-\de),u(t),u(t-\de)).
\]
With obvious modification, we can introduce notation with $b$ replaced by $\tilde{b},\;\si,\;\theta$ and $\ga$, and/or $x$ replaced by $x',\;u,\;u'$. For a sub-$\si$-field $\cG$, we use $\EE^\cG$ to denote the conditional expectation given $\cG$.

 We then define the Hamiltonian
\begin{eqnarray}\label{H}
&&\mathcal{H}(t, x, x',v, v', q,r,\bar{r},\alpha)\nonumber\\
&\doteq&    b(t,x,x',v,v' ) q +
 \sigma(t,x,x',v,v'  ) r \nonumber\\
 &&+\theta(t,x,x',v,v'  )   \bar{r}+l(t,x ,v,v'  )+h(t,x)\tilde{Q}
   \nonumber\\
 && + \int_{ \mathcal{E} }   \gamma( t,x,x',v,v',\zeta ) \alpha(\ze)  \nu(d \ze)  .
\end{eqnarray}

Finally, we formulate the adjoint equations as
\begin{eqnarray}\label{adjoint}
\left\{\!\!\begin{array}{rcl}
-dq(t) &= &\Big( b _x (t) q(t)+\sigma_x (t)r(t) +\theta_x (t)\bar{ r}(t) \\
&&+\int_{ \mathcal{E} }\gamma_x (t,\ze ) \alpha   (t,\zeta)\nu(d\zeta)+l_x (t)+ h_x(t)  \tilde{Q} (t) \\
&& +\bar{\mathbb E}^{\mathcal{F}_t}\Big[ b _{x'}(t+\de)q(t+\delta) +\sigma_{x'}(t+\delta)r(t+\delta) \\
&& +\theta_{x'}(t+\delta)\bar{ r}(t+\delta ) \nonumber \\
&& +\!\int_{ \mathcal{E} }\!\gamma_{x'}(t\!+\!\delta,\ze) \alpha(t\!+\!\delta,\zeta)\nu(d\zeta)\Big]  \Big) dt \\
&& -r(t)dB(t)-\bar{r}(t)dW^v(t)  -\int_{ \mathcal{E} }\alpha (t-,\zeta)\widetilde{N}(d\ze, dt), \\
 q(T) &= &\varphi_x(x(T)), ~ (q,r,\bar{r},\ga)(t)=0, ~ t\in (T,T+\delta],
\end{array}\right. \nonumber \\
\end{eqnarray}
and
\begin{equation}\label{P}
\left\{\begin{array}{rcl}
  -d P(t)&=&  \ell \(t, x(t),  v(t), v(t-\delta) \) dt-QdB(t) \\
  & & - \tilde{ Q}dW^v(t),  \nonumber\\
P(T)&=& \varphi(x(T)).
\end{array}\right.
\end{equation}

Here is our first main result.
\begin{theorem}\label{Necessary}
Assume that Hypotheses \ref{HP1} and \ref{HP2}  hold. Let $u(\cdot)$ be an  optimal
control and $  x(\cdot) $  the corresponding state.  Then, for any $v$ such that $v+u\in\cU_{ad}$, we have
\begin{equation}\label{smp}
\left\{\begin{array}{rcll}
\bar{\mathbb{E}}\left[\langle\mathcal{H}_{u}(t)
+ \mathcal{H}_{u'}(t+\delta),v(t)-u(t) \rangle \Big|\mathcal{G}_{t} \right]\le  0  
\quad \mbox{if }\quad 0\le t\le T-\de,\\
\bar{\mathbb{E}}\left[\langle\mathcal{H}_{u}(t),v(t)-u(t) \rangle
 \Big|\mathcal{G}_{t} \right]\le 0  \quad
 \mbox{if }\quad  T-\de<t\le T.
\end{array}\right.
\end{equation}
As a consequence, if $u(t)\in U^o$ (the interior of $U$), then
\begin{eqnarray}\label{smp1}
\left\{\begin{array}{rcll}
\bar{\mathbb{E}}\left[\mathcal{H}_{u}(t)
+ \mathcal{H}_{u'}(t+\delta) \Big|\mathcal{G}_{t} \right]
&=& 0&\quad\mbox{if }\quad 0\le t\le T-\de,\\
\bar{\mathbb{E}}\left[\mathcal{H}_{u}(t)
 \Big|\mathcal{G}_{t} \right]&=& 0&\quad \mbox{if } \quad T-\de<t\le T.
\end{array}\right. \nonumber \\
\end{eqnarray}
\end{theorem}

\begin{remark}
The above equations are identities involving the nonlinear filtering (or, equivalently, the conditional expectations) of the state and the adjoint processes. The study of their numerical approximation remains a challenging problem which we wish to come back in a future research. As a demonstration, we will provide an example in the last section where the state equation is linear and hence, the Kalman filtering is applicable. Actually, the separation principle also holds for that example.
\end{remark}

\section{Stochastic Maximum Principle}
\setcounter{equation}{0}
\renewcommand{\theequation}{3.\arabic{equation}}

In this section, we will present the proof of Theorem \ref{Necessary}.

For any $(\varepsilon,v')\in (0,1)\times
  \cU_{ad}$, let $v=v'-u$ and
 $   x^{u+\va v }(\cdot) $ be the solutions of
(\ref{subs-state}) with $v$ replaced by $ u+\va v  $.
Making use of Burkholder-Davis-Gundy inequality and Gronwall's inequality, we have
the following estimates (Lemmas \ref{lemma1}-\ref{lemma2}). 
\begin{lemma}\label{lemma1} Let Hypothesis \ref{HP1} hold. Then
\label{xZ}
\begin{eqnarray*}
\sup\limits_{0\leq t \leq T} \mathbb E |x^v(t)|^8 &\leq &C \bigg(1+\sup\limits_{-\delta\leq t\leq0}\mathbb E \( \xi(t)^8 +\eta(t)^8\)    + \sup\limits_{0 \leq t \leq T} \mathbb E\int_0^T |v(t)|^8 d t   \bigg),
\end{eqnarray*}
\[
\sup\limits_{0\leq t \leq T} \mathbb E |Z^v(t)|^ {   2}< \infty.
\]
\end{lemma}

\begin{lemma}\label{lemma2} Let Hypothesis \ref{HP1} hold. Then, there exits a constant $C$ such that
\label{xZ-per}
\[
\sup\limits_{0\leq t \leq T} \mathbb E |x^{u+\va  v}(t)-x^u(t)|^8\leq C \va ^8,
\]
and
\[
\sup\limits_{0\leq t \leq T}
  \mathbb E  |Z^{u+\va  v}(t)-Z^u(t)|^2  \leq C \va ^2.
\]
\end{lemma}

Since the proofs of Lemmas \ref{lemma1}-\ref{lemma2} are routine, we omit them.

Let
\[
\begin{array}{rl}
b^1(t)= & b_x(t) x^{1}(t)+b_{x'}(t) x^{1}(t-\delta)+b_u (t)v(t) \\
& +b_{u'}(t)v(t-\delta).
\end{array}
\]
$\si^1(t),\;\th^1(t),\;\ga^1(t,\ze)$ are defined similarly with $b$ replaced by $\si,\;\th,\;\ga(\cdot,\ze)$, respectively.

We introduce the following variational equations:
\begin{equation}\label{Z1}
\left\{
\begin{array}{rcl}
 dZ^{1}(t)&=&\Big[Z^{1}(t)h(t)+Z(t)h_x(t)x^{1}(t)  \Big]dY(t)\\
   Z^{1}(0)&=&\mbox{\ }0,
\end{array}
\right.
\end{equation}
and
\begin{equation}\label{x1}
\left\{\begin{array}{ccl}
 dx^{1}(t)&=&b^1(t) dt+ \si^1(t)dB(t)+\th^1(t)dW(t)\nonumber\\&&+\int_{ \mathcal{E }}\ga^1(t-,\ze)\widetilde{N}(d\zeta, dt) ,    \\
   x^{1}(0)&=&0,\;\; v(t)=0,\qquad\qquad t\in[-\delta, 0].
\end{array}
\right.
\end{equation}

For any $ v(\cdot)+u(\cdot) \in \mathcal{U}_{ad} $, it is easy to see that  under Hypothesis \ref{HP1}, (\ref{Z1}) and (\ref{x1}) admit a
unique solution, respectively. Furthermore,  we have the following lemma.

\begin{lemma}\label{limitxZ}
If Hypothesis \ref{HP1} holds, then
\begin{equation}\label{limit}
\begin{array}{rcl}
 & {\mathbb{E}}| x ^1(t)|^8 < \infty ,\quad
  {\mathbb{E}}|
  Z ^1(t)|^4 < \infty .\\
\end{array}
\end{equation}
\end{lemma}

\begin{lemma}\label{limitxyzxi}
If Hypothesis \ref{HP1} holds and
\[\tilde{\vartheta}^{\va}(t)=\frac{\vartheta^{u+\va v}(t)-\vartheta^{u}(t) }{\va }-\vartheta^{1}(t)\quad \hbox{with}\quad \vartheta=x,  Z.   \]
Then,
\begin{equation}\label{eq0506a}
\lim\limits_{\va\rightarrow0}\sup\limits_{0\leq t\leq T} {\mathbb{E}}|
 \tilde{x}^\va(t)|^4=0,
\end{equation}
\begin{equation}\label{Zvarepsilon}
\begin{array}{rcl}
 \lim\limits_{\va\rightarrow0}
 \sup\limits_{0\leq t\leq T}{\mathbb{E}}|
 \tilde{Z}^\va(t)|^2=0,
\end{array}
\end{equation}

\end{lemma}

Proof: Note that
\begin{equation}\label{eq0510a}
\begin{array}{rl}
d\tilde{x}^\va(t)= & \tilde{b}^{1,\va}(t)dt+\si^{1,\va}(t)dB(t)
+\th^{1,\va}(t)dY(t) \\
& +\int_{ \mathcal{E }}\ga^{1,\va}(t-,\ze)\widetilde{N}(d\zeta, dt),
\end{array}
\end{equation}
where
\[
b^{1,\va}(t)=\frac{1}{\va}\(b^{(\va)}(t)-b(t)\)-b^1(t),
\]
with
\[b^{(\va)}(t)=  b(t, x^{u+\va v} (t), x^{u+\va v}(t-\delta), u(t)+\va v(t),u(t-\delta)+\va v(t-\delta) ).\]
The other notation is defined by modifying $b^{(\va)}(t)$ in an obvious way.
By the boundedness of the second derivatives, we can then estimate
\[
\begin{array}{l}
|\tilde{b}^{1,\va}(t)| \le K\(|\tilde{x}^\va(t)|+|\tilde{x}^\va(t-\de)|\) \\
\quad\quad +\!K\va\(1\!+\!|x^{1}(t)|^2\!+\!|x^{1}(t\!-\!\de)|^2\!+\!|v(t)|^2\!+\!|v(t\!-\!\de)|^2\).
\end{array}
\]
The other coefficients can be estimated similarly. It then follows from (\ref{eq0510a}) and (\ref{limit}) that
 \begin{eqnarray}\label{esti-tidelx}
f(t)&\equiv&{\mathbb{E}}\sup_{s\le t}|\tilde{x}^\va(s)|^4\\
&\leq& K {\mathbb{E}}\int_0^t|\tilde{x}^\va(t)|^4 ds+ K {\mathbb{E}}\int_0^t|\tilde{x}^\va(t-\delta)|^4 ds+K\va^4.\nonumber
\end{eqnarray}
Combining (\ref{esti-tidelx}) with
\begin{eqnarray*}
&& {\mathbb{E}}\int_0^t|\tilde{x}^\va(s-\delta)|^4 ds
 ={\mathbb{E}} \int_{-\delta}^{t-\delta }|\tilde{x}^\va(s )|^4 ds      \\
&=& {\mathbb{E}} \int_{0}^{t-\delta }|\tilde{x}^\va(s )|^4 ds
\leq  {\mathbb{E}} \int_{0}^{t  }|\tilde{x}^\va(s )|^4 ds ,
\end{eqnarray*}
 we arrive at
\[f(t)\le K\int^t_0f(s)ds+\va^4K.\]
Then, Gronwall's inequality yields
\[
f(t)\le \va^4Ke^{Kt}\le K\va^4.
\]
This implies the desired results. 
\hfill$\Box$
\vspace{0.3cm}

Now we proceed to proving our main result.

{\em Proof of Theorem \ref{Necessary}:} As $u(\cdot)$ is optimal, we have
\begin{eqnarray}\label{dcostfunc}
0 &\geq & \frac{d}{d\va} J(\va)\mid_{  \va=0}\nonumber\\
&= & \lim \limits_{\va\rightarrow 0}\frac{ J(u(\cdot)+\va v(\cdot) )- J(u(\cdot))     }{\va}\nonumber
\\
 &=&{\mathbb{E}}\int_0^T\Big\{Z^{1}(t)l(t)+Z(t)\big[\ell_x(t)x^{1}(t)
+\ell_u(t)v(t)\nonumber\\
&& + \ell_{u'}(t)v(t-\delta)\big]\Big\}dt\nonumber\\
  &&+{\mathbb{E}}[Z(T)\varphi_x(x(T))x^{1}(T)]
 +{\mathbb{E}}[Z^{1}(T)\varphi(x(T))].\nonumber\\
&=& \bar{\mathbb E}\int_0^T\Big(\Gamma(t)\ell(t)+\ell_x(t)x^{1}(t)+\ell_u(t)v(t)\nonumber
\\&& +  \ell_{u'}(t)v(t-\delta)\Big)dt
 + \bar{\EE}  [\varphi_x(x(T))x^{1}(T)] \nonumber \\
&& +\bar{\EE}[\Gamma(T)\varphi(x(T))], 
\end{eqnarray}
where $ \Gamma (\cdot)=Z^1(\cdot)Z^{-1}(\cdot)$. By It\^{o}'s formula,
we have
\begin{eqnarray}\label{Gamma}
\left\{
\begin{array}{rcl}
  d\Gamma(t)&=&  h_x (t,x(t))x^1(t) \( dY(t)-h(t,x(t))dt \)\nonumber\\
 &= & h_x (t,x(t))x^1(t)dW^v(t), \nonumber \\
   \Gamma(0)&=&0.
\end{array}\right.
\end{eqnarray}

Applying It\^{o}'s formula and taking expectation, we obtain
\begin{eqnarray*}\label{x1qito}
&&\bar{\mathbb E}   [x^1(T) q (T)]\\
&=& \bar{ \mathbb E}   \int_0^T\Bigg(- \Big (\ell_x(t) +\tilde{Q}(t)h_x(t)\Big) x^1(t)
+\Big( b _u (t)v(t)+ b _{u'}(t)v(t-\delta)\Big) q(t) \nonumber\\&&+\Big(\sigma_u (t)v(t)+ \sigma_{u'}(t)  v(t-\delta) \Big)   r(t)  
+\Big(\theta_u (t)v(t)+ \theta_{u'}(t)  v(t-\delta) \Big)   \bar{r}(t)  
\nonumber\\ 
&&+ \int_{ \mathcal{E} }\Big(\gamma_u (t,\ze)v(t)+ \gamma_{u'} (t,\ze)v(t-\delta)\Big) \alpha(t,\ze)   \nu(d\zeta)\Bigg) dt\nonumber\\
&&+\int_0^T\Big( -x^1(t) \bar{\mathbb E}^{\mathcal{F}_t} [ b _{x'}(t+\de)q(t+\delta) ]  
+q(t) b_{x'}(t)x^1(t-\delta)                                                   \nonumber\\
&&    -x^1(t) \bar{\mathbb E}^{\mathcal{F}_t} [\sigma_{x'}(t+\delta)r(t+\delta)    ]+ r(t)\sigma_{x'}(t)x^1(t-\delta)   \nonumber\\
&&  -x^1(t) \bar{ \mathbb E}^{\mathcal{F}_t} [ \theta_{x'}(t+\delta)\bar{ r}(t+\delta )   ] + \bar{r}(t)\theta_{x'}(t)x^1(t-\delta)    \nonumber\\
&&  -x^1(t)\bar{\mathbb E}^{\mathcal{F}_t} [ \int_{ \mathcal{E} }\gamma_{x'}(t+\delta) \alpha   (t+\delta,\zeta)\nu(d\zeta)   ]        
+\alpha(t)\int_{ \mathcal{E}} \gamma_{x'}(t)x^1(t-\delta) \nu(d\zeta)\Big)dt .\nonumber
\end{eqnarray*}
Note that
\begin{eqnarray*}
&&\bar{\mathbb E} \int_0^T\!\!\( -x^1(t)\bar{ \mathrm E}^{\mathcal{F}_t} [ b _{x'}(t\!+\!\de)q(t\!+\!\delta)]\!+\!q(t) b_{x'}(t)x^1(t\!-\!\delta) \)dt                                                  \nonumber\\
&=&\bar{\mathbb E} \int_\delta^{T+\delta}\(-x^1(t-\delta)b_{x'}(t)q(t)  \)dt +\bar{\mathbb E} \int_0^T q(t) b_{x'}(t)x^1(t-\delta)dt \nonumber \\
&=&\bar {\mathbb E} \int_0^\delta b_{x'}(t)x^1(t-\delta)q(t)dt -\bar{\mathbb E} \int_\delta^{T+\delta}b_{x'}(t)x^1(t-\delta)q(t)dt=0,
\end{eqnarray*}
and other terms are similar. Hence, we can continue the above calculation with 
\begin{eqnarray*}
&&\bar{\mathbb E}   [x^1(T) q (T)]\\
 &=& \bar{ \mathbb E}   \int_0^T\Bigg(- \Big ( \ell_x(t) +\tilde{Q}(t)h_x(t)\Big) x^1(t)
+\Big( b _u (t)v(t)+ b _{u'}(t)v(t-\delta) \Big)q(t) \nonumber\\&&+\Big(\sigma_u (t)v(t)+ \sigma_{u'}(t)  v(t-\delta) \Big)   r(t)  
+\Big(\theta_u (t)v(t)+ \theta_{u'}(t)  v(t-\delta) \Big)   \bar{r}(t)  \nonumber\\
&&+ \int_{ \mathcal{E} }\Big(\gamma_u (t,\ze)v(t)+ \gamma_{u'} (t,\ze)v(t-\delta)\Big) \alpha(t,\ze)   \nu(d\zeta)\Bigg) dt.
\end{eqnarray*}
 By It\^o's formula and taking expectation again, we have
  \begin{eqnarray*}\label{Gammavarphito}
\bar{\mathbb E}   [\Gamma (T) \varphi\(x(T)\)] &=& - \bar{\mathbb E }  \int_0^T\(\Gamma(t) l(t)-\tilde{Q}   (t) h_x(t)x^{1}(t)\) dt.
 \end{eqnarray*}
Combining the above   equalities and (\ref{dcostfunc}),
   we have
 \[\bar{\mathbb E} \int_0^T \(\langle \mathcal{H}_u(t),v(t) \rangle
 +\langle  \mathcal{H}_{u'}(t), v(t-\delta) \rangle\)dt
 \leq 0. \]
 Namely,
  \begin{eqnarray*}
 \bar{\mathbb E}\int^T_{T-\de}\!\!\!\<\cH_u(s),v(s)\>ds +\bar{\mathbb E}\int^{T-\de}_0\!\!\!\<\cH_u(s)+\cH_{u'}(s+\de),v(s)\>ds\!
 \le 0,
 \end{eqnarray*}
 where we used the fact that for $s\in[-\de,0]$, \[v(s)=v'(s)-u(s)=\eta(s)-\eta(s)=0.\]
 
 For $t\in(T-\de,T)$, we can take
 \[v(s)=\left\{\begin{array}{ll}
 0&\mbox{ if }s\notin[t,t+\ep]\\
 v(s)&\mbox{ if }s\in[t,t+\ep],
 \end{array}\right.\]
 where $v-u\in \cU_{ad}$. Thus,
 \[\frac1{\ep}\bar{\mathbb E}\int^{t+\ep}_t\<\cH_u(s),v(s)\>ds\le 0.\]
 Taking $\ep\to 0$, we get $\bar{\mathbb E}\(\<\cH_u(t),v(t)\>|\cG_t\)\le 0$. The other conclusions can be proved similarly.
 \hfill $\Box$
\vspace{0.3cm}

In some cases,  the observation functions 
are linear and hence not bounded. From a careful check of the proofs of Theorem \ref{Necessary}, we have

\begin{remark}\label{remark1} \rm
The conclusions remain true when $h$ being bounded in Hypothesis \ref{HP1} is replaced by
\[
\mathbb E  \( | Z(t) |^2 +| Z^1(t) |^4 \)  <\infty.
\]
However, this condition will not cover general linear model which is being studied in another research using method specifically using the linearity.  If the diffusion coefficients in the state equation  are constants and the admmisible controls are restricted to linear feedback ones, this condition can be verified.
\end{remark}

\section{A sufficient maximum principle}
\setcounter{equation}{0}
\renewcommand{\theequation}{4.\arabic{equation}}

Inspired by \cite{lenhartxiongyong2016}, we now establish a sufficient condition without assumption of the concavity of $\mathcal{H} $. The result is derived based on the expansion of the cost functional associated with  (\ref{state})
and (\ref{observation}).

\begin{theorem}\label{Sufficiency} 
For $ u (\cdot)\in \mathcal{U}_{ad} $, and let  $v(\cdot)$ be such that $u(\cdot) +\varepsilon  v(\cdot)\in \mathcal{U}_{ad} $.
The following Taylor expansion holds:
\begin{eqnarray}\label{Jtalor12}
   J(u\!+\!\varepsilon v  )=J(u)\!+\!\varepsilon J_1(u;v )\!+\! \varepsilon^2 J_2 (u;v )\!+\! o(\varepsilon^2),
\end{eqnarray}
where
\begin{eqnarray*}\label{J1}
  J_1(u;v )  &=& \bar{\mathbb E}\int_0^T\Big(\Gamma(t) \ell(t)+\ell^1(t)\Big)dt
 \\
 &&+ \bar{\mathbb E}  [\varphi_x(x(T))x^{1}(T)]
 +\bar{\mathbb E}[\Gamma(T)\varphi(x(T))],
\end{eqnarray*}
with
\[\ell^1(t)=\ell_x(t)x^{1}(t)
+\ell_u(t)v(t)\nonumber\\
 + \ell_{u'}(t)v(t-\delta),\]
and
\begin{eqnarray*}\label{J2}
&&J_2(u;v )  =\bar{\mathbb E}\int_0^T\(\Gamma^1(t)\ell(t)+\Gamma(t)\ell^1(t)+  \ell_x(t)x^2(t)+\ell^2(t) \)    dt \\
&& \quad\quad\quad\quad +   \bar{\mathbb E} \big [\varphi_x(x(T))x^{2}(T)+\varphi_{xx} (x(T))\(x^{1}(T) \)^2\nonumber\\
&& \quad\quad\quad\quad + \Gamma(T)\varphi_x(x(T))x^1(T)  + \Gamma^1(T)\varphi(x(T))  \big],
\end{eqnarray*}
with
\begin{eqnarray}\label{Gamma-1}
  d\Gamma^1(t)&=&   \( h_{xx} (t,x(t))\(x^1(t)\)^2 \!\!+\!h_{x} (t,x(t)) x^2(t) \!\)dW^v(t), \nonumber\\
   \Gamma^1(0)&=&0,
\end{eqnarray}
 \begin{equation}\label{x2-suf}
 \left\{\begin{array}{ccl}
 dx^{2}(t)&=&b^2(t)dt+\si^2(t)dB(t)+\th^2(t)dW(t)\nonumber\\&&+\int_{\cE}\ga^2(t-,\ze)
 \widetilde{N}(d\zeta, dt) ,    \\
   x^{2}(0)&=&0, v(t)=0,\quad t\in[-\delta, 0],
\end{array}\right.
\end{equation}
and
\begin{eqnarray*}
 b^2(t)&=&
b_{xx}(t)x^1(t)+ b_{xx'} (t)x^1(t-\delta)+b_{xu}(t)v(t)\nonumber\\
&&+b_ {xu'}(t)v(t-\delta)     \big)x^{1}(t)+b_x(t)x^2(t) \\
 &&+
  \big(b _{x'x}(t)x^1(t) + b _{x'x'} (t)x^1(t-\delta)  +b _{x'u}(t)v(t)\nonumber\\&& +b _ {x'u'}(t)v(t-\delta)     \big)x^{1}(t-\delta)+b_{x'} (t)x^2(t-\delta) \\
  &&
  +\big(b_{ux}(t)x^1(t)+ b_{ux'} (t)x^1(t-\delta)+b_{uu}(t)v(t) \\&&+b_ {uu'}(t)  v(t-\delta)   \big)v(t)\\
 &&+
  \big(b _{u'x}(t)x^1(t)+ b _{u'x'} (t)x^1(t-\delta)+b _{u'u}(t)v(t)\\&&+b _ {u'u'}(t) v(t-\delta)    \big)v(t-\delta).
  \end{eqnarray*}
The other coefficients are defined similarly.

As a consequence, if
\[  
J_2\( u;v \)\le 0, \quad \forall  v(\cdot)\in \(\mathcal{U}_{ad}-u\)\setminus\{0\},   
\]
and
\[ 
J_1\( u;v \)= 0.
\]
Then $ u(\cdot)$ is a locally optimal control.
\end{theorem}
Proof: Similar to Lemma \ref{limitxyzxi}, we can prove that
\[x^{u+\ep v}(t)=x(t)+\ep x^1(t)+\ep^2 x^2(t)+o(\ep^2).\]
Further expansion of $\phi=b,\ \si,\ \th,\ \ga$, we obtain that
\[\phi^\ep(t)=\phi(t)+\ep\phi^1(t)+\ep^2\phi^2(t)+o(\ep^2).\]
Finally, we can take Taylor expansion for $J(u+\varepsilon v)$ at $\ep=0$ to arrive at the desired result.
\hfill $\Box$

\section{Two examples}
\setcounter{equation}{0}
\renewcommand{\theequation}{\thesection.\arabic{equation}}

In this section we present two examples to
demonstrate the applicability of Theorem \ref{Necessary}. In the first example, the optimal strategy is provided. In the second one concerns the optimal investment problem with delay. Finally, the numerical solution is shown.

\subsection{Example 1}
 
 Consider the following problem with state equation:
\begin{eqnarray}\label{state1}
\left\{\begin{array}{rcl}
dx (t)& =& ( \sin(x (t))+ a x(t-\delta)+ bv(t) )    dt \\
 && + \sigma  x (t) dB(t)  + \theta dW^v(t) + \!\int_{\mathcal{E}}\!\gamma\widetilde{N}(d\zeta, dt), \hfill t\in[0,T],   \nonumber  \\
x (t) &=& \xi(t), \quad v(t)=\eta(t),\quad\quad  t\in[-\delta,0],
\end{array}\right.  \end{eqnarray}
where    $a, \ b, \ \sigma,\ \th,  \gamma $ are constant.
The observation process $Y(t)$ is given by
\begin{eqnarray*}
\left\{\begin{array}{rcl}
dY(t) & \!\!=\!\! & \cos( x (t))dt+dW^v(t),\nonumber\\
Y(0)  & \!\!=\!\! & 0.\nonumber
\end{array}\right.
\end{eqnarray*}
The performance functional is assumed to have the form
\begin{eqnarray*}
\begin{array}{l}
J(v(\cdot)) = \bar {\mathbb{E}}\Big[\int _0^T v^2(t) dt+x(T)\Big].\\
\end{array}  \nonumber \\
\end{eqnarray*}
Then, we have 
\begin{eqnarray*}
&&\mathcal{H}(t, x, x',v,  q,r,\bar{r},\alpha)\nonumber\\
&\doteq &  ( \sin(x (t))+ a x(t-\delta)+ bv(t)) q +
 \sigma x(t) r \nonumber\\
 &&+\theta   \bar{r} +v^2(t)+\cos(x(t))\tilde{Q}
  +   \int_{ \mathcal{E} }   \gamma  \alpha(\ze)  \nu(d \ze).
\end{eqnarray*}
In addition, (\ref{adjoint}) reads
\begin{eqnarray}\label{qexam1}
 \left\{\!\!\begin{array}{rcl}
-dq(t) &= &\Big( \cos(x (t) )q(t)+\sigma r(t)  -\sin(x(t) ) \tilde{Q} (t) \\
&&  +\mathbb E^{\mathcal{F}_t}\Big[a q(t+\delta) \Big]  \Big) dt \\
&& -r(t)dB(t)-\bar{r}(t)dW(t)  -\int_{ \mathcal{E} }\alpha (t-,\zeta)\widetilde{N}(d\ze, dt), \\
 q(T) &=& 1, ~ (q,r,\bar{r},\ga)(t)=0, ~ t\in (T,T+\delta],
\end{array}\right. \nonumber \\
\end{eqnarray}
and (\ref{P}) becomes
\begin{eqnarray} \label{Pexam2}
\left\{\begin{array}{rcl}
  -d P(t)&=& v^2(t)  -QdB(t)  - \tilde{ Q}dW^v(t),  \nonumber \\
P(T)&=&  x(T) . \nonumber
\end{array}\right.
\end{eqnarray}
According to Theorem \ref{Necessary},   for any $v$ such that $v+u\in\cU_{ad}$, we have
\[  
\bar{\mathbb{E}}\left[\langle  b q(t)+2u(t),v(t)-u(t) \rangle \Big|\mathcal{G}_{t} \right]\le 0. 
\]
 If $u(t)\in U^o$ (the interior of $U$), then
\[  \bar{\mathbb{E}}\left[\langle  b q(t)+2u(t) ,v(t)-u(t) \rangle \Big|\mathcal{G}_{t} \right]=  0.           \]
By the arbitrariness of $v(t)$, we obtain the optimal control
\[ 
u(t)= -    \bar{\mathbb{E}} [ q(t) \Big|\mathcal{G}_{t}],  
\]
where $q$ satisfies (\ref{qexam1}). Thus, we adopt our stochastic maximum principle to convert the nonlinear control problem into the calculation of filter which can be calculated using filtering techniques.

\subsubsection{Example 2}

 Suppose  a financial market consisting of a risk-free asset (bond) and a risky asset (stock). Specifically, the price process of the risk-free asset is given by
\begin{eqnarray*}
d S_{0}(t)  &=&   r_0(t)  S_{0}(t)d t,\quad r_0 > 0,
\end{eqnarray*}
where $r_0(t)$ is the published interest rate
and the price process of the risky assets follows the following stochastic differential equation \begin{eqnarray*}\label{surplus}
d S_{1}(t)  &=&  \mu(t)  S_{1}(t) d t   +  \sigma
 (t) S_{1}(t) dW (t),
\end{eqnarray*}
 where  $ \mu(t)\;\ (>r_0(t) )$ is the appreciation rate process,  which is not observable directly. $\sigma  (t) $ is the volatility coefficient, and $W(t) $ is a
standard Brownian motion.

A strategy $v(t)$  represents the amount invested in the  risky asset at time $t$.
Naturally, it is assumed to  be adapted to an observation generated sub-filtration $\cG_t$.

Without delay, the  controlled wealth process is given by
\begin{eqnarray*}\label{eq0521a}
d x(t) &=&\(r_0(t)\(x(t)-v(t)\)\) dt + \mu(t)v(t)dt\nonumber\\
&&+ \sigma(t) v(t) dW (t).
\nonumber
\end{eqnarray*}

In real situations, there are some delays between the investment behavior and change of the wealth process. Thus, our model should be modified as:
\begin{eqnarray*}\label{eq0521a}
d x(t) &=&\(r_0(t)\(x(t-\delta)-v(t)\)\) dt + \mu(t)v(t)dt\nonumber\\
&&+ \sigma(t) v(t) dW (t).
\nonumber
\end{eqnarray*}
 The policymaker can get information from the stock price
\begin{eqnarray*}
\left\{\begin{array}{rcl} dY (t)&=& \frac{1}{ \sigma(t )} \(  \mu (t) -\frac{1}{2}  \sigma ^2(t)  \) dt+dW(t),\\
Y(0)&=& 0,
\end{array}
\right.
\end {eqnarray*}
where $Y(t) =\log S(t)/\sigma(t)$.

Similar to Theorem 3.1 in \cite{zx}, using separation principle, we have
\begin{eqnarray*}\label{X1}
\left\{\begin{array}{rcl} dx (t)
&=& [r_0(t)\(x(t\!-\!\delta)\!-\!v(t)\)] dt \!+\! \hat{\mu}(t)v(t)dt \\
&  & + \sigma(t) v(t) d \nu(t),\nonumber \\
x(0)&=& x_0,
\end{array}
\right.
\end{eqnarray*}
where $ \hat{\mu}(t) =\mathrm E[ \mu(t)| \mathcal{G}_t ]     $, and the innovation process
\[
d\nu(t)= \frac 1{\sigma(t)}d \log S(t) - \frac{1}{\sigma(t)} \(\hat{\mu}(t) -\frac{1}{2}\sigma^2 (t)  \)  dt.
\]
The objective is to maximize
\begin{eqnarray}
J(v(\cdot))=\frac{1}{2}\mathbb{E} \left[\int_{0}^{T} -( v (t)-a(t) )^2 dt+  x
 (T)  \right].\nonumber
\end{eqnarray}
This means that the investor should not only prevent $v(t)$ from large deviation with respect to a certain level $a(t)$, but also maximize the terminal wealth.

In this setting, from  (\ref{H})  and  (\ref{adjoint}), the Hamiltonian becomes
\begin{eqnarray*}\label{H-examp}
\begin{array}{l}
\mathcal{H}(t, x, x',v, v', q,r,\bar{r})  \\
= (r_0(t)\( x(t-\delta) -v(t)) +\hat{\mu}(t) v(t) \)   q(t) \\
\quad + \sigma(t) v(t) r(t) - \frac{1}{2}\( v(t)-a(t)  \)^2 ,
\end{array}
\end{eqnarray*}
and the adjoint equation is given by
\begin{equation}\label{adjoint-exam}
\left\{
\begin{array}{rcl}
-d q(t)&=& \mathbb E^{\mathcal{F}_t} \( r_0 (t)   q(t+\delta)\) dt  - r (t)d\nu(t)    \nonumber\\ q(T) & =& \frac{1}{2}.
\end{array}
\right.
\end{equation}
Applying the stochastic maximum principle (\ref{smp}),
we get the optimal control satisfying
\begin{equation}\label{eq1204a}
  v(t)= \( \hat{\mu}(t)-r_0(t) \)  q(t)+\sigma (t)  r (t)  -a(t),
\end{equation}
where $(q,r)$ is the unique solution to BSDE (\ref{adjoint-exam}).
Applying filtering technique to the adjoint equation (\ref{adjoint-exam}), we arrive at the BSDE (\ref{adjoint-exam}) for $(\hat{q},\hat{r})$.

As a further demonstration of the applicability of our results, we consider the numerical solution of (\ref{surplus}) with $\sigma (t)= 0.3\sin(2t)+0.1$, and
 \[d \mu(t) = \alpha \mu(t) dt +\beta d \bar W(t),\]
 where  $\alpha $ and $\beta $ are constant,  and $\bar W(t) $ is a Brownian motion independent of $W$.
 According to the Kalman-Bucy filtering theory,
\begin{equation}\label {filtermu}
d\hat \mu(t)=\alpha \hat\mu (t)dt+\frac{\gamma(t)}{\sigma(t)}d\bar W(t),   \end{equation}
where $\gamma(t)$ is the solution to the following Riccati equation
\begin{equation}\label{Reccati}
\frac{d\gamma(t)}{dt}=2\alpha \gamma(t)-\frac{\gamma^2(t)}{\sigma^2(t) }+\beta^2.\end{equation}

The terminal time $T$ is set to $T=1$, and the volatility rate of the stock price  $\sigma (t) $ is set to  $ \sigma(t) =0.3 \sin(2 t )+0.1$ and $a(t)=0.1t$.
We now  give a discretization technique to deal with  numerical scheme for the optimal control problem, which is along the same line with \cite{zhang2019}.
Let
\begin{eqnarray*}\label{partition}
 |\pi|  &:=&\frac{\delta}{m},  \nonumber\\
 \pi  &:=& t_0<t_1=\frac{\delta}{m}< t_2=\frac{2\delta}{m}\cdots< t_{n-1}=\frac{(n-1)\delta}{m} \\
 & & \leq t_{n}=T
 \end{eqnarray*}
 be a partition of $[0,T]$,
Define $\pi(t)\triangleq   t_{i-1}$, for $ t \in [t_{i-1}, t_i)$.  In order to compute (\ref{adjoint-exam}), let
\begin{eqnarray*}\label{D-y}
q_{t_n}^\pi  &=& q(T), \nonumber \\
   q_{t}^\pi  &=&  q_{t_i}^\pi+\(  \mathbb E ^{\mathcal{F}_{t_i}}  q_{t_i+\delta}^\pi       \)  \Delta t_i-\int_t^ {t_i} r_s^\pi d \nu_s,\nonumber\\
    &&
 \quad
     t\in [t_{i-1},t_i),i=n,n-1,...,1.
\end{eqnarray*}
Then, we focus our attention to the integral term. Using Lemma 2.7 in \cite{zhang2004},  we set
 \begin{eqnarray*}
   r_{i-1}^\pi   =  \frac{1}{\Delta t_{i}}
\mathbb E \( q_i^\pi  \( \nu_{t_i}- \nu_{t_{i-1}} \)      \)
.\nonumber
 \end{eqnarray*}

Firstly, we use Matlab to solve (\ref{Reccati}) and get the path of $\gamma$. Then,  discretizing equation (\ref{filtermu}) using Euler scheme, we obtain an approximation of the filtering process  $\hat{\mu}_t$. That is, \begin{eqnarray*}\label{disc-hatmu}
   \hat\mu_{t_i}^{\pi}  &=& \hat\mu_{t_{i-1}}^\pi+ \alpha \hat\mu ^\pi_{t_{i-1}}\Delta t_i+ \frac{\gamma(t_{i-1})}{\sigma(t_{i-1})} \(\bar W_{t_i}-\bar W _{t_{i-1}}\). \nonumber \end{eqnarray*}
In the following, computing the value of $q^\pi_{t_n}$. At last, recall (\ref{eq1204a}), we have the path of the control $v$.

The plots for  $\gamma(t)$ and $\hat{\mu}(t)$  are shown in the following figure.
\begin{figure}[tbh!]
\centering
\hspace{-1.8cm} \quad\quad
\includegraphics[width=2.97in]{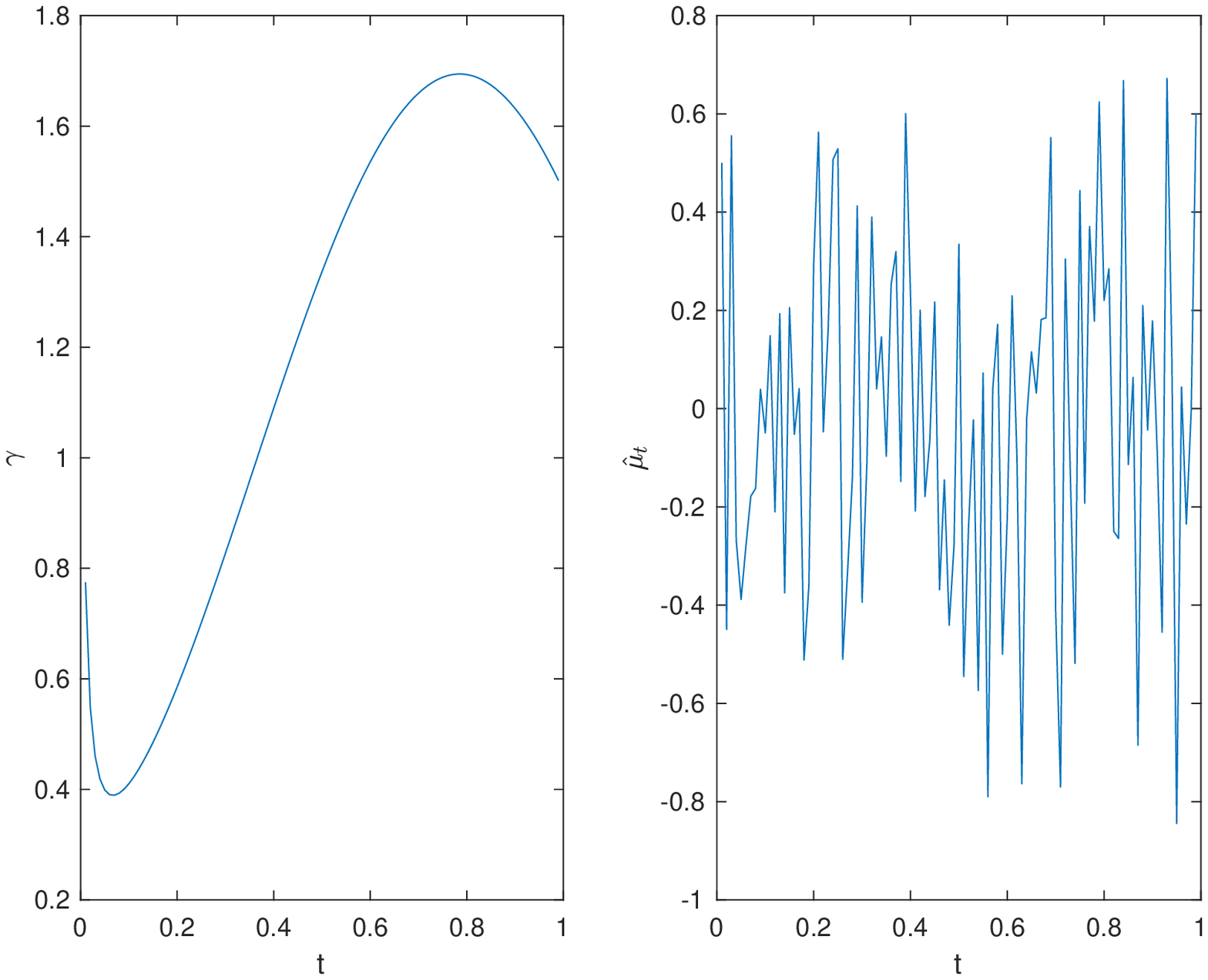} \\
Figure 1. 
\caption{The path\label{plot_beta}}
\end{figure}

In order to obtain the effect of delay,  comparison of  value function is made. Since the valuation function is a expectation,  the simulations are repeated $100$ times.

\begin{table}[tbh!]
\begin{center}
\begin{tabular}[c]{ |c| c| c| c|}
\hline
 &   $ \mathbb{E}  - \int _0^T (v(t)-a(t))^2 dt $  & $\mathbb{E} x(t)$  &    $ J(v)$ \\ \hline
  $\delta=0.48$  &-4.5873  & 2.0301 & -2.5573 \\ \hline
  $\delta=0.44$  &-4.3558  &  2.0900 & -2.2658 \\ \hline
  $\delta=0.42$ &-3.7601  &  2.0890 & -1.6711 \\ \hline
  $\delta=0.40$ &-3.1938  &  2.0874 & -1.1063 \\ \hline
  $\delta=0.38$ &-2.7869  &  2.0874 & -0.6919 \\ \hline
\end{tabular} \\
\caption{The effect of the delay term }
\label{table1}
\end{center}
\end{table}
One point of interest is the performance of the deviation from   $v(t)$ and psychological expectation level $a(t) $
 becomes smaller as delay becomes smaller, and the value function  becomes larger. In short, the delay leads to decreased value function. 

\section*{Acknowledgment}

The authors would like to thank two anonymous referees for constructive suggestions which substantially improved the presentation of the paper.

\end{document}